\documentclass[12pt, reqno]{amsart}
\usepackage{amsmath, amsthm, amscd, amsfonts, amssymb, graphicx, color}
\usepackage[bookmarksnumbered, colorlinks, plainpages]{hyperref}
\hypersetup{colorlinks=true,linkcolor=red, anchorcolor=green, citecolor=blue, urlcolor=red, filecolor=magenta, pdftoolbar=true}

\renewcommand{\phi}{\varphi}
\newtheorem{theorem}{Theorem}[section]
\newtheorem{lemma}[theorem]{Lemma}

\newtheorem{proposition}[theorem]{Proposition}

\theoremstyle{definition}
\newtheorem{definition}[theorem]{Definition}

\theoremstyle{remark}
\newtheorem{remark}[theorem]{Remark}

\numberwithin{equation}{section}

\def\AA{{\mathcal A}}

\def\GG{{\mathcal G}}

\def\MM{{\mathcal M}}

\def\TT{{\mathcal T}}

\def\bbR{\mathbb{R}}
\def\bbC{\mathbb{C}}
\def\bbZ{\mathbb{Z}}

\def\bbH{\mathbb{H}}

\begin{document}
	
	\setcounter{page}{1}
	
	\title[  Algebras of Toeplitz Matrices with Quaternion Entries]{Algebras of Toeplitz Matrices with Quaternion Entries}
	
	\author[ Muhammad Ahsan Khan Ameur Yagoub and]{Muhammad Ahsan Khan $^1$, Ameur Yagoub $^2$}
	\address{$^{1}$  Department of Mathematics, University of Kotli Azad Jammu $\&$ Kashmir, Kotli 11100, Azad Jammu $\&$ Kashmir, Pakistan.}
	\email{\textcolor[rgb]{0.00,0.00,0.84}{ahsankhan388@hotmail.com}}

	\address{$^{2}$ Laboratoire de math\'ematiques pures et appliqu\'es. Universit\'e de Amar Telidji. Laghouat, 03000. Algeria.}
	\email{\textcolor[rgb]{0.00,0.00,0.84}{a.yagoub@lagh-univ.dz}}
	\thanks{The first author is the corresponding author, ahsankhan388@hotmail.com.}
	
	
	
	
	\subjclass[2010]{15B33, 15B05, 11R52 }
	
	\keywords{Quaternions; Toeplitz matrices; maximal algebras; circulant matrices}
	
	
	\begin{abstract}
	The classification of maximal left algebras of quaternion Toeplitz matrices is a harder problem that has received little attention up to now. In this paper, we introduce certain families of maximal left algebras of Toeplitz matrices with entries from an algebra of quaternions that cover various classes of the left algebras of quaternion Toeplitz matrices.
	\end{abstract} \maketitle
\section{Introduction}
A square matrix is called Toeplitz if every descending diagonal from left to right is constant. In other words, the elements of the matrix are arranged in such a way that each row is a shifted version of the previous row.

	These matrices arise naturally in many areas of Mathematics and are important both in theory and applications. For instance, it is well known that a large class of matrices is similar to Toeplitz matrices \cite{Heinig2001, Mackey1999}. Moreover, it is proved that every matrix can be expressed as a product of Toeplitz matrices (see \cite{Ye2016}). Apart from this, these matrices have some of the most attractive computational properties and are amenable to a wide range of disparate algorithms. We refer the reader to \cite{Widom1965} for a detailed study of these matrices.   
	
	Quaternions are a fascinating Mathematical construct and extend the concept of complex numbers. Instead of two components like complex numbers, quaternions have four components. These numbers play important roles across many areas of Mathematics generally as algebraic systems, signal processing, differential geometry, and quantum mechanics, etc.
	
	Matrices over commutative rings take attention but, matrices having noncommuting entries (quaternion entries) have not been investigated very much yet. This is basically due to intrinsic algebraic difficulties that appear with respect to their non-commutativity. During the last decade, a large amount of research has been concentrated on Toeplitz matrices over the field of complex numbers, while their study over quaternions is quite negligible.

 Hamilton first introduced the set of real quaternions (see \cite{Hamilton1953}). While the seminal work concerning commutative quaternions was first presented by \cite{Segre1892}. Kosal and Tosun \cite{Kosal2014, Kosal2015} investigated some algebraic properties of commutative quaternion matrices using complex representations of commutative quaternion matrices. We refer the reader to \cite{Kosal2014, Kuipers, Cayley, Hamilton, Kleyn2014, Shuster1, Fallon, Shuster, Kosal2015} and \cite{Zhang1997, Altmann, Tait} for a detailed study of quaternions and their matrices.  In \cite{Gu2003}, the authors provide basic properties of Toeplitz and Hankel matrices over the algebra of complex numbers; most of the results therein deal with the products of these structured matrices, which, in general, are not structured over the algebra of quaternions. The most usual and basic reference for complex Toeplitz matrices is Grenander and Szeg¨o \cite{Grenander1958}.

	The collection of quaternion Toeplitz matrices is not closed with respect to the multiplication of matrices. So it is interesting to find classes of quaternion Toeplitz matrices that have the structure of left vector space as well as the ring. The general problem of characterizing the left algebras (or right algebras) of quaternion Toeplitz matrices is a very hard problem and no work has been done hitherto. The purpose of the present paper is to obtain the classification of maximal left algebras of quaternion Toeplitz matrices. It is probably too tough to hope for a complete classification, but the purpose is to identify possible classes of such left algebras for quaternion Toeplitz matrices. 
	
	This paper is structured as follows: 
	By means of Section 2, we want to make sure that
	the reader has become acquainted with quaternions and their algebraic properties required when we start the main work in upcoming sections. In Section 3 we will introduce Toeplitz matrices over an algebra of quaternions.  In the last Section, we introduce a certain class of left maximal algebras of quaternion Toeplitz matrices and prove some fundamental results concerning it.
	
	 \section{Quaternions and their Basic Algebraic Properties}
	In this section, the main object of study is the set of real quaternions and their algebraic properties. Hamilton introduced the Hamiltonian quaternions for representing vectors in space, \cite{Hamilton1953, Hamilton1969}.
	\begin{definition}
		The set of real quaternions denoted by $\mathbb{H}$ and is defined as 
		\[
		\mathbb{H}:=\Big\{\alpha=\alpha_0+\alpha_1i+\alpha_2j+\alpha_3k: \alpha_0,\alpha_1,\alpha_2,\alpha_3\in\bbR\Big\}
		\]
		where $i,j,k\notin\bbR$ are versors satisfying the following multiplication rule: 
		$$i^2=j^2=k^2=ijk=-1$$.
	\end{definition}
	The name $\bbH$ is given in honor of an Irish Mathematician William Rowan Hamilton who first discovered these in the 1840s.  
	Since quaternion arithmetic is defined by the behavior of versors, We can also derive lots of other relations from these, for instance, $ij=k$, $ki=j$, $jk=i$, $ji=-k$, $ik=-j$, $kj=-i$ and $jki=-1$. 
	As an additive group $\bbH$ is isomorphic to four copies of $\bbR$, i.e., $\bbR\oplus \bbR\oplus \bbR\oplus \bbR$. The map $\alpha_0+\alpha_1i+\alpha_2j+\alpha_3k\longrightarrow (\alpha_0, \alpha_1, \alpha_2, \alpha_3)$ is clearly a group isomorphism of $\bbH$ onto $\bbR\oplus \bbR\oplus \bbR\oplus \bbR$.
	
	The real quaternions are the obvious generalization of complex numbers, We define the addition  and multiplication of quaternions as follows:\\
	If $\alpha=\alpha_0+\alpha_1i+\alpha_2j+\alpha_3k$,  $\beta=\beta_0+\beta_1i+\beta_2j+\beta_3k$ are  $\mathbb{H}$, then 
	$$\alpha+\beta=(\alpha_0+\beta_0)+(\alpha_1+\beta_1)i+(\alpha_2+\beta_2)j+(\alpha_3+\beta_3)k$$ and 
	\begin{align*}
		\alpha\beta&=(\alpha_0\beta_0-\alpha_1\beta_1-\alpha_2\beta_2-\alpha_3\beta_3)\\
		&+(\alpha_0\beta_1+\alpha_1\beta_0+\alpha_2\beta_3-\alpha_3\beta_2)i\\
		&+(\alpha_0\beta_2+\alpha_2\beta_0+\alpha_1\beta_3-\alpha_3\beta_1)j\\
		&+(\alpha_0\beta_3+\alpha_3\beta_0+\alpha_1\beta_2-\alpha_2\beta_1)k.
	\end{align*}
	It is easy to see that with respect to the above operations, $\bbH$ is a division ring and also that under multiplication quaternions are not commutative. Due to this reason, one must take some care in order to perform the multiplication of quaternions. Under the usual operation of addition and scalar multiplication, $\bbH$ is a four-dimensional vector space over $\bbR$. 
	
	The following theorem from \cite{Zhang1997} summarizes some of the main algebraic properties of quaternions. 
	\begin{theorem}
		Let $\alpha=\alpha_0+\alpha_1i+\alpha_2j+\alpha_3k$ and $\beta=\beta_0+\beta_1i+\beta_2j+\beta_3k $ in $\bbH$, then  \begin{itemize} 
			\item [(i)] Every $\alpha$  can  be expressed in a unique way as  $\alpha=\gamma_0+\gamma_1j$, $\gamma_0,\gamma_1\in\bbC$; 
			\item[(ii)] In general, $(\alpha+\beta)^{2} \neq \alpha^{2}+2\alpha\beta+\beta^{2}$;
			\item[(iii)] $\alpha^2+1=0$ has  infinitely  many  roots over $\bbH$.
		\end{itemize}
		
	\end{theorem}
	It is also notable that one can express any quaternion in terms of a $2\times 2$ matrix having complex entries. Let $\alpha=\alpha_0+\alpha_1i+\alpha_2j+\alpha_3k$ is in $\bbH$. Expressing the versors $1, i,j$, and $k$ as 
	\[
	1=\begin{pmatrix}
		1 & 0\\
		0 & 1
	\end{pmatrix}, \quad i=\begin{pmatrix}
		i & 0\\
		0 & -i
	\end{pmatrix},\quad j=\begin{pmatrix}
		0 & 1\\
		-1 & 0
	\end{pmatrix}, \quad k=\begin{pmatrix}
		0 & i\\
		i & 0
	\end{pmatrix}.
	\]
	Then simple computation imply that $i^2=j^2=k^2=-1$ and $ij=k$, $jk=i$, and $ki=j$. So $\alpha$ takes the following form 
	\begin{align*}
		\alpha&=\alpha_0+\alpha_1j+\alpha_2j+\alpha_3k\\
		&=\begin{pmatrix}
			\alpha_0+\alpha_1i & \alpha_2+\alpha_3i\\
			-\alpha_2+\alpha_3i & \alpha_0-\alpha_1i
		\end{pmatrix}.
	\end{align*}
 We now quote from \cite{Altun2021}, the definition of the inner product on $\bbH^n$. 
\begin{definition}
Let $n\in\bbZ^+$, and  $x=\begin{pmatrix}
    x_0\\
    x_1\\
    \vdots\\
    x_{n-1}
\end{pmatrix}$ and $y=\begin{pmatrix}
    y_0\\
    y_1\\
    \vdots\\
    y_{n-1}
\end{pmatrix} $ be in $\bbH^n$. The quaternion valued function $<\cdot,\cdot>:\bbH^n\times\bbH^n\longrightarrow \bbH $ given by  $$<x,y>=\sum_{k=0}^{n-1}\overline{y_k}\ x_k$$
defined an inner product on $\bbH^n$.
\end{definition}
With this inner product $\bbH^n$ is a left inner product space over $\bbH$ and likewise the unitary space $\bbC^n$, the set of vectors $$e_0=\begin{pmatrix}
	1\\
	0\\
	\vdots\\
	0
\end{pmatrix}, e_1=
\begin{pmatrix}
	0\\
	1\\
	\vdots\\
	0
\end{pmatrix},\cdots,  e_{n-1}=
\begin{pmatrix}
	0\\
	0\\
	\vdots\\
	1
\end{pmatrix}$$
form an orthonormal basis for $\bbH^n$(these vectors also serve as the standard Hamel Basis for $\bbH^n$).
	
	\section{Toeplitz Matrices with Quaternion Entries}
	Throughout we label the indices of any matrix $A$ from $0$ to $n-1$. 
	We denote by $\MM_n \left[\bbH\right]$, the set of all square matrices with quaternion entries. Thus if $A\in\MM_{n}[\bbH]$ then $A=(\alpha_{rs})_{r,s=0}^{n-1}$, with $\alpha_{rs}\in\bbH$ for every $0\leq r,s\leq n-1$.
	
	If $A=(\alpha_{rs})_{r,s=0}^{n-1}$, and, $B=(\beta_{rs})_{r,s=0}^{n-1}$ are in $\MM_n\left[\bbH\right]$ and $\alpha\in\bbH$, then define matrix addition and scalar multiplication component-wise as follows:
	\begin{align*}
		A+B&=(\alpha_{rs}+\beta_{rs})_{r,s=0}^{n-1},\\
		\alpha A&=(\alpha\alpha_{rs})_{r,s=0}^{n-1}.
	\end{align*}
	With respect to the above defining operations, $\MM_n \left[\bbH\right]$ is a left vector space over $\bbH$. If one multiply $A=(\alpha_{rs})_{r,s=0}^{n-1}\in\MM_n[\bbH]$ by $\alpha\in\bbH$ from right then, $\MM_n\left[\bbH\right]$ is also a  right vector space over $\bbH$. Now we define a Toeplitz matrix, whose entries all belong to $\bbH$.
	
	A finite square matrix is called a quaternion Toeplitz matrix if its entries along each negative sloping diagonal are constant. That is, the matrix $A=(\alpha_{rs})_{r,s=0}^{n-1}$ is quaternion Toeplitz if $\alpha_{r_1,s_1}=\alpha_{r_2,s_2}$ whenever $r_1-s_1=r_2-s_2$, for all $r_1,s_1,r_2,s_2=1-n, 2-n, \cdots, -1,0,1,\cdots, n-1$. 
	
	From the definition, it is clear that a quaternion Toeplitz matrix of size $n^2$ depends upon $2n-1$ parameters $1-n, 2-n, \cdots -1,0,1,\cdots,n-1.$ The word ``quaternion” refers to the fact that in the above matrix representation, the entries are from the algebra of quaternions. Thus if $A\in\MM_n[\bbH]$ is Toeplitz then it has the following structure:
	\[
	A=
	\begin{pmatrix}
		\alpha_0& \alpha_{-1}& \alpha_{-2} & \cdots &\alpha_{1-n}\\
		\alpha_1& \alpha_0 & \alpha_{-1} & \cdots &\alpha_{2-n}\\
		\alpha_{2}& \alpha_{1} & \alpha_{0} & \cdots &\alpha_{3-n}\\
		\vdots & \vdots & \vdots &\ddots & \vdots\\
		\alpha_{n-1} & \alpha_{n-2}& \alpha_{n-3} & \cdots &\alpha_{0}
	\end{pmatrix}, \quad \alpha_{r}\in\bbH \quad\hbox{for all} \quad 1-n \leq r\leq n-1.
	\]
	We denote by  $\TT_{n}[\bbH]$ the set of all Toeplitz matrices with entries from $\bbH$. If $A=(\alpha_{r-s})_{r,s=0}^{n-1}$ is in $\TT_{n}[\bbH]$ then $A$ is called quaternion  circulant matrix if $\alpha_{r}=\alpha_{r-n}$ for every $1\leq r\leq n-1$.

 Let $S$ be the matrix consisting of zeros except for ones below the principal diagonal, i.e, 

\[S=\begin{pmatrix}
 		0 & 0&0 &\cdots &0&0\\
 		1 & 0&0 & \cdots & 0&0\\
        0&1&0&\cdots &0&0\\
 		\vdots & \vdots & \vdots&\ddots & \vdots\\
 		0 &0&0& \cdots &1&0
 	\end{pmatrix}
  \]

  Then its its adjoint $S^{\ast}$ is the matrix given as 

\[
S^{\ast}=\begin{pmatrix}
 		0 & 0&0 &\cdots &0&0\\
 		1 & 0&0 & \cdots & 0&0\\
        0&1&0&\cdots &0&0\\
 		\vdots & \vdots & \vdots&\ddots & \vdots\\
 		0 &0&0& \cdots &1&0
 	\end{pmatrix}
  \]
It is  clear that both $S$ and $S^*$ are in $\TT_n[\bbH]$. 
Recall that if  $x$ and $y$ are in $\bbH^n$, then the tensor product $x\otimes y$ is a square matrix of size $n^2$ over $\bbH$ and is defined as  $x\otimes y(z)=<z,y>x, \quad \hbox{for every}\quad x,y\in\bbH^n$.

The following result characterized all the matrices of $\TT_n[\bbH]$ among the matrices of $\MM_n[\bbH]$. We are adding its detailed proof just for completeness as it is similar to the proof of Lemma 2.2 of \cite{Gu2003}.
\begin{proposition}
\label{Toeplitz classificatiois siliar n}
 $A\in\MM_n[\bbH]$ is in $\TT_{n}[\bbH]$ if and only if there exist vectors $x$ and $y$ in $\bbH^{n}$ such that 
\[
A-SA S^*=x\otimes e_{0}+e_{0}\otimes y.
\]
\end{proposition}
\begin{proof}
Suppose that $A$ is in $\TT_n[\bbH]$, then $A$ has the form $A=(\alpha_{r-s})_{r,s}^{n-1}$. Simple computation yields that 
\[
A-SA S^*=
 \begin{pmatrix}
\alpha_{0}  &   \alpha_{1}  &   \alpha_{2}   &\cdots &  \alpha_{n-1}\\
\alpha_{-1}&    0  &    0    &     \cdots &     0\\
\alpha_{-2}&    0  &0   &    \cdots      &0\\
\vdots &    \vdots & \vdots      &\ddots & \vdots\\
 \alpha_{1-n} &     0     &0      & \cdots &0
 \end{pmatrix}. 
\]

 If one take $x=\begin{pmatrix}
 			\alpha_{0}\\
 			\alpha_{-1}\\
            \vdots\\
 		     \alpha_{1-n}
\end{pmatrix}$ and $ y=\begin{pmatrix}
 			0\\
 			\alpha_{1}\\
 		    \vdots\\
            \alpha_{n-1}
\end{pmatrix}$, then it is easy to see that 
\[
A-SA S^*=x\otimes e_{0}+e_{0}\otimes y.
\]

In order to establish converse, let $A=(\alpha_{rs})_{r, s=0}^{n-1}$ be any arbitrary $n\times n$ matrix such that for some 
$x=\begin{pmatrix}
 			\alpha_{0}\\
 			\alpha_{-1}\\
            \vdots\\
 		     \alpha_{1-n}
\end{pmatrix}$ and $ y=\begin{pmatrix}
 			0\\
 			\alpha_{1}\\
 		    \vdots\\
            \alpha_{n-1}
\end{pmatrix}$  in $\bbH^n$, the identity  
$
A-SA S^*=x\otimes e_{0}+e_{0}\otimes y
$
holds. 
In standard basis of $\bbH^n$, one has

\begin{align*}
x\otimes e_{0}(e_{r})=
\begin{cases}\displaystyle\sum_{k=0}^{n-1}\alpha_{k}e_{k},\quad & \mbox{for}\quad r=0\\ 
0 \quad& \mbox{for}\quad  r=1,2,\cdots,n-1,
\end{cases}.
\end{align*}
and  
\[
e_{0}\otimes y(e_{r})=\delta_{r}e_{0}, \quad \hbox{for}\quad r= 0,1,\cdots,n-1.
\]

Therefore 
\begin{equation}\label{x}
x\otimes e_{0}+e_{0}\otimes y=\begin{pmatrix}
 		\alpha_{0}+\delta_{0}&\delta_{1}&\cdots &\delta_{n-1}\\
 		\alpha_{1}& 0 & \cdots & 0\\
 		\vdots & \vdots &\ddots & \vdots\\
 		\alpha_{n-1} &0& \cdots &0
 	\end{pmatrix}.
  \end{equation}
On the other hand one note that 
\begin{equation}\label{y}
 A-SA S^*=\begin{pmatrix}
 		\alpha_{00} & \alpha_{01}&\cdots &\alpha_{0,n-1}\\
 		\alpha_{10}&\alpha_{11}& \cdots &\alpha_{1,n-1}-\alpha_{0,n-2}\\
 		\vdots & \vdots &\ddots & \vdots\\
 		\alpha_{1-n,0} &\alpha_{n-1,1}-\alpha_{n-2,0}& \cdots &\alpha_{n-1,n-1}-\alpha_{n-2,n-2}
 	\end{pmatrix}.
\end{equation}
Comparing corresponding entries of (\ref{x}) and (\ref{y}) , one sees that 
$\alpha_{r,s}=\alpha_{r-1,s-1}$ for every  $1\leq r,s\leq n-1$ . This shows that $A$ is in $\TT_n[\bbH]$, which is what we wanted to prove.
\end{proof}
	\begin{proposition}\label{Toeplitz subspace}
		$\TT_{n}[\bbH]$ is a  subspace of a left vector space $\MM_{n}[\bbH]$.
	\end{proposition}
	\begin{proof}
		The proof is immediate by noting that if $A=(\alpha_{r-s})_{r,s=0}^{n-1}$, $B=(\beta_{r-s})_{r,s=0}^{n-1}$ and $\gamma\in\bbH$, then 	
		$
		A+\gamma B=(\alpha_{r-s}+\gamma\beta_{r-s})_{r,s=0}^{n-1}$ is in $\TT_{n}[\bbH]$.
	\end{proof}
	Note that $A=(\alpha_{r-s})_{r,s=0}^{n-1}$ with
	\begin{equation*}
		\alpha_{r-s}=\begin{cases}
			i\quad \hbox{if}\quad r-s=1,\\
			j\quad \hbox{if} \quad r-s=-1,\\
			0 \quad \hbox{otherwise}, 
		\end{cases}
	\end{equation*}
	is purely a quaternion Toeplitz the matrix, but the product 
	\[A^2=\begin{pmatrix}
		k& 0 & 0 & 0 & \cdots & 0\\
		0 & 0 & 0 & -1 & \cdots  & 0\\
		-1 & 0 & 0 & 0 & \cdots & 0\\
		\vdots& \vdots & \vdots & \vdots& \cdots &\vdots\\
		0 & 0 & 0 & 0  &\cdots & -k
	\end{pmatrix}\]is not a quaternion Toeplitz matrix. Thus likewise the case of complex entries, the product of two quaternion Toeplitz matrices needs not be a quaternion Toeplitz matrix. The following proposition gives us the precise criteria in this direction. 
	\begin{proposition}\label{le:basic condition for product toeplitz}
		Suppose that $A=(\alpha_{r-s})_{r,s=0}^{n-1}$
		and $B = (\beta_{r-s})_{r.s=0}^{n-1}$ are Toeplitz matrices with quaternion entries then $AB$ is in $\TT_{n}[\bbH]$ if and only if \begin{equation}\label{eq:basic product condition}
			\alpha_r\beta_{s-n}=\alpha_{r-n}\beta_s\quad\hbox{for all}\quad r,s=1,2,\dots n-1.
		\end{equation}
	\end{proposition}
	\begin{proof}
		Suppose that $AB$ is in $\TT_n[\bbH]$. Let us denote the product $AB$ by $C=(\gamma_{r,s})_{r,s=0}^{n-1}$, then for every $r,s=1,2,\dots ,n-1$, we have 
		\begin{equation}\label{sum1}
			\gamma_{r,n-s}
			=\displaystyle\sum_{k=0}^{n-1}\alpha_{r-k}\beta_{k+s-n}
		\end{equation}
		\begin{equation}\label{sum2}
			\gamma_{r-1,n-s-1}=
			\displaystyle\sum_{k=0}^{n-1}\alpha_{r-k-1}\beta_{k+s-n+1}	
		\end{equation} 
		Subtracting \eqref{sum1} and \eqref{sum2} yields
		\begin{equation}\label{sum3}
			\gamma_{r,n-s}-\gamma_{r-1,n-s-1}=\alpha_r\beta_{s-n}-\alpha_{r-n}\beta_r
		\end{equation}
		Since the product $AB$ is a quaternion Toeplitz matrix, then its elements along the negative sloping diagonals have the same value, consequently from \eqref{sum3} 
		\[ 
		\alpha_r\beta_{s-n}
		=\alpha_{r-n}\beta_r \quad \hbox{for every}\quad r,s=1,2,\dots ,n-1.
		\]
		Conversely suppose that the identity \eqref{eq:basic product condition} is true, then for every $r,s=1,2,\dots, n-1$, we have 
		\begin{align*}
			\gamma_{r,n-s}-\gamma_{r-1,n-s-1}
			&=\alpha_r\beta_{s-n}-\alpha_{r-n}\beta_r \\
			&=0.
		\end{align*}
		This shows that $AB$ is a  quaternion Toeplitz matrix.
	\end{proof}
 \begin{proposition}
		If $A=(\alpha_{r-s})_{r,s=0}^{n-1}$ and $B=(\beta_{r-s})_{r,s=0}^{n-1}$ are in $\TT_n[\bbH]$  with commuting entries then  $A$ and $B$ commute with each other.
	\end{proposition}
	\begin{proof}
		We have for every $r,s=0,1,\cdots, n-1,$
		\begin{align*}
			(AB)_{r,s}
			&= \sum_{k=0}^{n-1}\alpha_{r-k}\beta_{k-s}\\
			(BA)_{r,s}
			&= \sum_{k=0}^{n-1}\beta_{r-k}\alpha_{k-s}.
		\end{align*}
		Since entries are commuting, then rewriting the sum by denoting $k^\prime=r+s-k$, we obtain 
		\[
		(BA)_{r,s}= \sum_{k=0}^{n-1}\alpha_{k-s}\beta_{r-k}=\sum_{k^\prime=r+s-(n-1)}^{r+s}\alpha_{r-k^\prime}\beta_{k^\prime-s}
		\]
		If $r+s=n-1$, then the above sum is the same as the formula for $(AB)_{r,s}$. So $(BA)_{r,s}=(AB)_{r,s}$. Now, for instance,  let us suppose that $r+s<n-1$. Then only part of the sum is the same, and from the rest, we obtain
		\begin{align*}
			(AB)_{r,s}-(BA)_{r,s}
			&= \sum_{k=r+s+1}^{n-1}\alpha_{r-k}\beta_{k-s}- \sum_{r+s-(n-1)}^{-1}\alpha_{r-k}\beta_{k-s}\\
			&=\sum_{k=r+s+1}^{n-1}(\alpha_{r-k}\beta_{k-s}-\alpha_{r-k+n}\beta_{k-s-n})
		\end{align*}
		Applying Proposition \ref{eq:basic product condition}, we get $(AB)_{r,s}-(BA)_{r,s}=0$. The proof for the case $r+s>n-1$ can be established in a similar fashion.
  \end{proof}
  The converse of the above proposition is not true in general. This is because if one takes $$A=\begin{pmatrix}
      0 & 0 & 0 &\cdots &1\\
      0 & 0 & 0 &\cdots &0 \\
      0 & 0 & 0 & \cdots & 0\\
      \vdots & \vdots & \vdots& \ddots & \vdots\\
      1 & 0& 0 & \cdots &0 
  \end{pmatrix}
   \quad\hbox{and}\quad B=\begin{pmatrix}
      0 & 0 & 0 &\cdots &k\\
      0 & 0 & 0 &\cdots &0 \\
      0 & 0 & 0 & \cdots & 0\\
      \vdots & \vdots & \vdots& \ddots & \vdots\\
      k & 0& 0 & \cdots &0 
  \end{pmatrix}.$$	
  Then $A$ and $B$ are commuting quaternion Toeplitz matrices but their product 
  \[
AB=\begin{pmatrix}
      k & 0 & 0 &\cdots &0\\
      0 & 0 & 0 &\cdots &0 \\
      0 & 0 & 0 & \cdots & 0\\
      \vdots & \vdots & \vdots& \ddots & \vdots\\
      0 & 0& 0 & \cdots &k
  \end{pmatrix}.
  \]
  is not a quaternion Toeplitz matrix. 
  \section{Algebras of Toeplitz Matrices with Quaternions Entries}
	Fixing an algebra $\AA$ in $\bbH$ and $p,q\in\AA^\prime$. where $\AA^\prime$ denotes the commutant of $\AA$ (the set of all quaternions commuting with each element of $\AA$). It is easy to see that $\AA^\prime$ is also an algebra.  We symbolize by $\TT_{n}\left[\AA\right]$ the left vector space of Toeplitz matrices contained in $\TT_n[\bbH]$. We define the family by
	\begin{equation}\label{definition of algebra}
		\GG_{p,q}[\AA]:=\Big\{A=(\alpha_{r-s})_{r,s=0}^{n-1}| \alpha_r\in\AA,\quad p\alpha_{r-n}=q\alpha_{r}\Big\}
	\end{equation}
	We will use the following simple Lemma.  
	\begin{lemma}\label{kernel}
		Suppose that $p$, $q\in\bbH$ be fixed. If $\alpha$ is any arbitrary element of $\bbH$ such that $p\alpha=q\alpha=0$ , then $\alpha=0.$ 
	\end{lemma}
	\begin{proof}
		Suppose that $p$, $q\in\bbH$, then we may write 
		$$p=\begin{pmatrix}
			p_0+p_1i & p_2+p_3i\\
			-p_2+p_3i & p_0-p_1i
		\end{pmatrix}$$ and 
		$$q=\begin{pmatrix}
			q_0+q_1i & q_2+q_3i\\
			-q_2+q_3i & q_0-q_1i
		\end{pmatrix}.$$ Let $\alpha=\begin{pmatrix}
			\alpha_0+\alpha_1i & \alpha_2+\alpha_3i\\
			-\alpha_2+\alpha_3i & \alpha_0-\alpha_1i
		\end{pmatrix}$ be any arbitrary element of $\bbH$, such that $p\alpha=q\alpha=0$. Then if at least one of $p$ or $q$ is nonzero then  $\alpha=0$.  
	\end{proof}
	The above Lemma enables us to prove the following main result of this section. 
	\begin{theorem}\label{algebras}
		The family $\GG_{p,q}[\AA]$ form a left algebra in $\TT_{n}[\AA]$. 
	\end{theorem}
	\begin{proof}
		A simple straightforward verification shows that $\GG_{p,q}[\AA]$ is a left subspace of $\TT_{n}[\bbH]$. We need to only show that it is closed up to the usual multiplication of matrices. For this, let us suppose that $A$ and $B$ be any two arbitrary elements of $\GG_{p,q}[\AA]$, then we must have 
		\begin{equation}\label{eq:GpqforTU}
			\begin{split}
				& A=(\alpha_{r-s}), \quad p\alpha_{r-n}=q\alpha_{r}\\
				& B=(\beta_{r-s}),\quad  p\beta_{r-n}=q\beta_{r}.\\
			\end{split}
		\end{equation}  
		We have \begin{align*}
			(AB)_{r,s}-(AB)_{r+1,s+1}&=\sum_{k=0}^{n-1}\alpha_{r-k}\beta_{k-s}\\
			&=\alpha_{r-n+1}\beta_{n-1-s}+\alpha_{r+1}\beta_{-1-s}.\\
		\end{align*}
		Using formulas \ref{eq:GpqforTU} and multiplying with $p$ and $q$, it follows that 
		\begin{align*}
			p\left[(AB)_{r,s}-(AB)_{r+1,s+1}\right]&=0\\
			q\left[(AB)_{r,s}-(AB)_{r+1,s+1}\right]&=0
		\end{align*}
		By applying Lemma \ref{kernel} we obtained that 
		\[
		(AB)_{r,s}-(AB)_{r+1,s+1}=0. 
		\]
		Consequently  $AB$ is in $\TT_{n}[\AA]$. We denote the product $AB$ by $C=(\gamma_{r,s})_{r,s=0}^{n-1}$, and since $p,q\in\AA^\prime$, we have 
		\[
		q\gamma_{r}=q(AB)_{r,0}=q\sum_{k=0}^{n-1}\alpha_{k-r} \beta_k= \sum_{k=0}^{n-1}\alpha_{k-r} p\beta_{k-n}=p\sum_{k=0}^{n-1}\alpha_{k-r} \beta_{k-n}=p\gamma_{r-n}.
		\]  
		Therefore the product $AB$ is in $\GG_{p,q}[\AA]$. This proves, along with the previous fact that $\GG_{p,q}[\AA]$ 
		is a left vector space, that $\GG_{p,q}[\AA]$ is a left algebra in $\TT_{n}[\AA]$.
	\end{proof}
	\begin{remark}
		One sees that the algebra $\GG_{p,q}[\AA]$ covers various classes of the algebras of $\TT_{n}[\AA]$, as is seen below:
		\begin{itemize}
			\item If $p=q=1$, then $\GG_{p,q}[\AA]$ is the left algebra of all  quaternion circulant matrices. 
			\item If $p=0$, then $\GG_{0,q}[\AA]$ is the left algebra of all upper triangular quaternion Toeplitz matrices.
			\item Similarly, if $q=0$, then $\GG_{p,0}[\AA]$ is the left algebra of all lower triangular quaternion Toeplitz matrices.
			\item  $\GG_{p,q}[\AA]$ obviously contains the algebra of diagonal quaternion Toeplitz matrices.
		\end{itemize} 
	\end{remark}
	\begin{proposition}
		If $\AA$ is commutative then $\GG_{p,q}[\AA]$ is commutative.  
	\end{proposition}

	The following result describes when two left algebras of this type are equal. 
	\begin{theorem}
		Suppose that there also exist $\breve{p}, \breve{q}\in\AA$, then $\GG_{p,q}[\AA]=\GG_{\breve{p},\breve{q}}[\AA]$ if and only if $p\breve{q}=\breve{p} q$.
	\end{theorem} 
	\begin{proof}
		Suppose that $p\breve{q}=\breve{p} q$, and $A=(\alpha_{r-s})_{r,s=0}^{n-1}$ is in $\GG_{p,q}[\AA]$, then
		$ p\alpha_{r-n}=q\alpha_{r}$. Multiplying the identity $p\alpha_{r-n}=q\alpha_{r}$ with $\breve{p}$, we have 
		\begin{align*}
			\breve{p} p\alpha_{r-n}&=\breve{p} q\alpha_{r}\\
			&= p\breve{q} \alpha_{r}
		\end{align*} 
		or \[
		p(\breve{p} \alpha_{r-n}-q\alpha_r)=0.
		\] 
		In the same way, multiplying both sides of $p\alpha_{r-n}=q\alpha_{r}$ with $\breve{q}$, we get
		\[
		q(\breve{p} \alpha_{r-n}-q\alpha_r)=0. 
		\]
		Lemma \ref{kernel} then imply that $\breve{p} \alpha_{r-n}=\breve{q}\alpha_r$.  Consequently $A\in\GG_{\breve{p},\breve{q}}[\AA]$ and therefore $\GG_{p,q}[\AA]\subset \GG_{\breve{p},\breve{q}}[\AA]$. Now for the reverse inclusion, let $A=(\alpha_{r-s})_{r,s=0}^{n-1}\in\GG_{\breve{p},\breve{q}}[\AA]$, then we have  $\breve{p}\alpha_{r-n}=\breve{q}\alpha_{r}$. Multiplying the equation $\breve{p}\alpha_{r-n}=\breve{q}\alpha_{r}$  by $q$ and using $p\breve{q}=\breve{p} q$,  we get $p \alpha_{r-n}=q\alpha_r$. Therefore, $\GG_{\breve{p},\breve{q}}[\AA]\subset \GG_{p,q}[\AA]$.
		
		Conversely suppose that $\GG_{p,q}[\AA]= \GG_{\breve{p},\breve{q}}[\AA]$. We will show that $p\breve{q}=\breve{p} q$. Since the matrix 	
		$$B=\begin{pmatrix}
			0 & q & q\cdots &q\\
			p & 0& q\cdots &q\\
			p & p &0 \cdots &q\\
			\vdots& \vdots & \vdots\ddots& \vdots  \\
			p& p & p \cdots &0 
			
		\end{pmatrix}$$ is in $\GG_{p,q}[\AA]=\GG_{\breve{p},\breve{q}}[\AA]$, and therefore $p\breve{q}=\breve{p} q$. Finishing the proof.
	\end{proof}
	
	The next result is related to the maximality of $\GG_{p,q}[\AA]$ as a left algebra in $\TT_{n}[\AA]$.
	\begin{theorem}
		$\GG_{p,q}[\AA]$ is maximal left algebra in $\TT_n[\AA]$ if and only if $\left\{p,q\right\}^\prime=\AA$. 
	\end{theorem} 
	\begin{proof}
		Suppose that $\GG_{p,q}[\AA]$ is  maximal left algebra in $\TT_{n}[\AA]$. We will show that $\left\{p,q \right\}^\prime=\AA$. Let $\left\{p,q \right\}^\prime=\breve{\AA}$. Since $p,q$ are in the commutant of $\AA$ then by definition $\left\{p,q\right\}^\prime \supset \AA$ and as a consequence, we have $\GG_{p,q}[\AA]\subset\GG_{p,q}[\breve{\AA}]$. This shows that $\GG_{p,q}[\AA]$ is not a maximal left algebra, which is a contradiction to our assumption. Thus the maximality of $\GG_{p,q}[\AA]$ imply that $\left\{p,q\right\}^\prime=\AA$. Conversely suppose that $\left\{p,q\right\}^\prime=\AA$ and $\GG$ be any arbitrary left algebra of quaternion Toeplitz matrices such that $\GG_{p,q}[\AA]\subseteq \GG$. Let $G=(g_{r-s})_{r,s=0}^{n-1}$
		be any arbitrary element of $\GG$. Since $\left\{p,q\right\}^\prime=\AA$, then the matrix
		$$\Tilde{G_{p,q}}=\begin{pmatrix}
			0 & q & q\cdots &q\\
			p & 0& q\cdots &q\\
			p & p &0 \cdots &q\\
			\vdots& \vdots & \vdots\ddots& \vdots  \\
			p& p & p \cdots &0 
		\end{pmatrix}$$ is in $\GG_{p,q}[\AA]$. Since $\GG_{p,q}[\AA]\subseteq \GG$ then the product $\Tilde{G_{p,q}} G$
		is in $\GG$. Since $\GG$ is a left algebra of quaternion Toeplitz matrices, comparing entries along the diagonals yields that $pg_{r-n}=qg_r$. This shows that $G\in\GG_{p,q}[\AA]$ and hence $\GG=\GG_{p,q}[\AA]$, i.e., $G_{p,q}[\AA]$ is maximal in $\TT_2[\AA]$. This is what we wanted to prove. 
	\end{proof}
	\section*{Declarations}
	The authors confirm that there are no known conflicts of interest associated with this manuscript and there has been no significant financial support for this work that could have influenced its outcome. 
	\section*{Data Availability} This manuscript has no associated data.
	\section*{Acknowledgments} 

	\bibliographystyle{amsplain}

\end{document}